\documentclass[12pt]{amsart}
\usepackage{amsfonts, amsmath, amssymb, amscd, latexsym, graphicx}

\newtheorem{thm}{Theorem}[section]

\newtheorem{lem}[thm]{Lemma}

\newtheorem{defn}[thm]{Definition}

\def\qed{\hfill\rlap{$\sqcup$}$\sqcap$\par}

\begin{document}

\title{Metric properties of higher-dimensional Thompson's groups}

\author{Jos\'e Burillo}
\address{
Departament de Matem\`atica Aplicada IV, Escola Polit{\`e}cnica
Superior de Castelldefels, Universitat Polit\`ecnica de Catalunya,
08860 Castelldefels, Barcelona, Spain} \email{burillo@mat.upc.es}

\author{Sean Cleary}
\address{Department of Mathematics,
The City College of New York \& The CUNY Graduate Center, New York,
NY 10031} \email{cleary@sci.ccny.cuny.edu}

\thanks{The authors are grateful for the hospitality of the
Centre de Recerca Matem\`atica. The first author acknowledges
support from MEC grant \#MTM2006-13544-C02-02 and is grateful for
the hospitality of The City College of New York.
This work was supported (in part) by a grant from The City University of New York PSC-CUNY Research Award Program.}

\date{10 October 2008}

\begin{abstract}
Higher-dimensional Thompson's groups $nV$ are finitely presented groups
described by Brin which generalize dyadic self-maps of the unit interval to
dyadic self-maps
of $n$-dimensional unit cubes. We describe some of the metric properties of
higher-dimensional Thompson's groups.  We give descriptions of elements
based upon tree-pair diagrams and  give upper and lower bounds for word length in
terms of the size of the diagrams.  Though these upper and lower bounds are
somewhat separated, we show that there are elements realizing the lower
bounds and that the fraction of elements which are close to the upper bound
converges to 1, showing that the bounds are optimal and that the upper bound
is generically achieved.
\end{abstract}

\maketitle

\section*{Introduction}

Thompson's groups serve as a wide range of interesting examples
of unusual group-theoretic behavior.  The family of Thompson's
groups includes the original groups described by Thompson
(commonly denoted $F$, $T$ and $V$) as well as generalizations
in many different directions.  The bulk of these generalizations includes
groups which can be regarded as self-maps of the unit interval.
  Brin \cite{brinhigher1,brinhigher2}
describes higher-dimensional groups $nV$, generalizations of
Thompson's group $V$, which are described naturally in terms
of dyadic self-maps of $n$-dimensional cubes.  Little is known about
these groups aside from their simplicity.  Brin describes these
elements geometrically in terms of dyadic interpolations of collections
of dyadic blocks, and gives presentations for the group $2V$.

We describe elements in higher-dimensional Thompson's groups
as being given by tree-pair diagrams.  Though tree-pair diagrams
usually take advantage of the natural left-to-right ordering of subintervals
of the unit interval, by using several types of carets, Brin gives a
a description of  elements of $nV$ via tree-pair diagrams.  A natural question
is how the size of these tree-pair diagrams corresponds to the word
length of elements with respect to finite generating sets.

We give upper and lower bounds for the word lengths of elements of
higher-dimensional Thompson's groups $nV$ respect to finite
generating sets, in terms of the size of tree-pair descriptions of
elements.   An element with a reduced tree-pair description of size $N$ has
word length between $\log N$ and $N \log N$ (up to the standard
affine equivalences) with respect to the
standard finite generating set.  This of course also thus holds
for any finite generating set.

The authors are grateful for helpful conversations with Matt Brin,
Melanie Stein and Claire Wladis.

\section{Background on higher-dimensional Thompson's groups}

 Brin \cite{brinhigher1} describes  higher-dimensional Thompson's groups $nV$,
 giving presentations and showing that the group $2V$ is simple.
 Here we summarize the relevant properties of the the groups $nV$ for
 use in this article.

We denote the unit interval $[0,1]$ as $I$, and correspondingly
$I^n$ denotes the $n$-dimensional unit cube.

An element of $V$ is given by two finite rooted binary
trees with the same number of leaves and a permutation, which represents a bijection between the
leaves of the two trees. Hence an element of $V$ can be seen as a
triple $(T_+,\pi,T_-)$, where $T_+$ and $T_-$ are trees with $k$
leaves, and $\pi\in \mathcal S_k$.

Each binary tree can be seen as a way of subdividing the interval
$I$ into dyadic subintervals of the type $[\frac{i}{2^r},\frac{i+1}{2^r}]$,
where $r>0$ and $0\le i < 2^r$.
A rooted binary tree gives instructions for successive
halvings of subintervals to obtain a particular dyadic subdivision.
Given two binary trees with $n$ leaves, and a permutation in $\mathcal S_n$,
an element in $V$ is represented as
a left-continuous map of the interval into itself, sending each
interval in the first subdivision to a corresponding interval
in the second subdivision,  as specified by the given permutation.

For more details on the group $V$, including presentations and a
proof of its simplicity, see Cannon, Floyd and Parry \cite{cfp} as well as Brin \cite{brinhigher1}.

To obtain an element in $2V$ we will define a partition of $I^2$ into dyadic
rectangles of the type
$$
\left[\frac{i}{2^r},\frac{i+1}{2^r}\right]\times
\left[\frac{j}{2^s},\frac{j+1}{2^s}\right].
$$
We can again regard such partitions as being obtained by a
successive halving process, but now there are two possible
halvings that can occur at each stage to a specified dyadic rectangle---
a horizontal or a vertical subdivision.
We can obtain a refinement of any dyadic partition of $I^2$
via iterated horizontal and vertical subdivisions.  This
gives means of describing elements of $2V$ as
pairs of dyadic subdivisions, each with $m$ rectangles,
and a permutation in $S_m$ giving the bijection between
the rectangles.

We can represent elements of
$2V$ also with pairs of binary trees. Since there are two types of
subdivisions,  vertical and horizontal, we will
consider binary trees which contain two types of carets, ``vertical"
and ``horizontal" carets, represented by ``triangular" and ``square"
carets, respectively.

\begin{figure}
\includegraphics[width=3in]{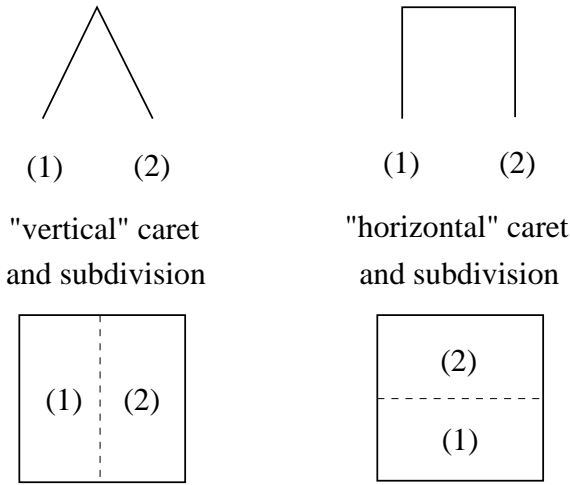}\\
\caption{The two types of carets and their corresponding
subdivisions} \label{carets}
\end{figure}

To fix a convention, we will represent a vertical subdivision with a
traditional triangular caret, where the left and right leaves
naturally represent the left and right rectangles. A horizontal
subdivision will be represented by a square caret, and in it, the
left leaf represents the bottom rectangle, and the right leaf
represents the top rectangle, as shown in Figure \ref{carets}.

\begin{figure}
\includegraphics[width=6in]{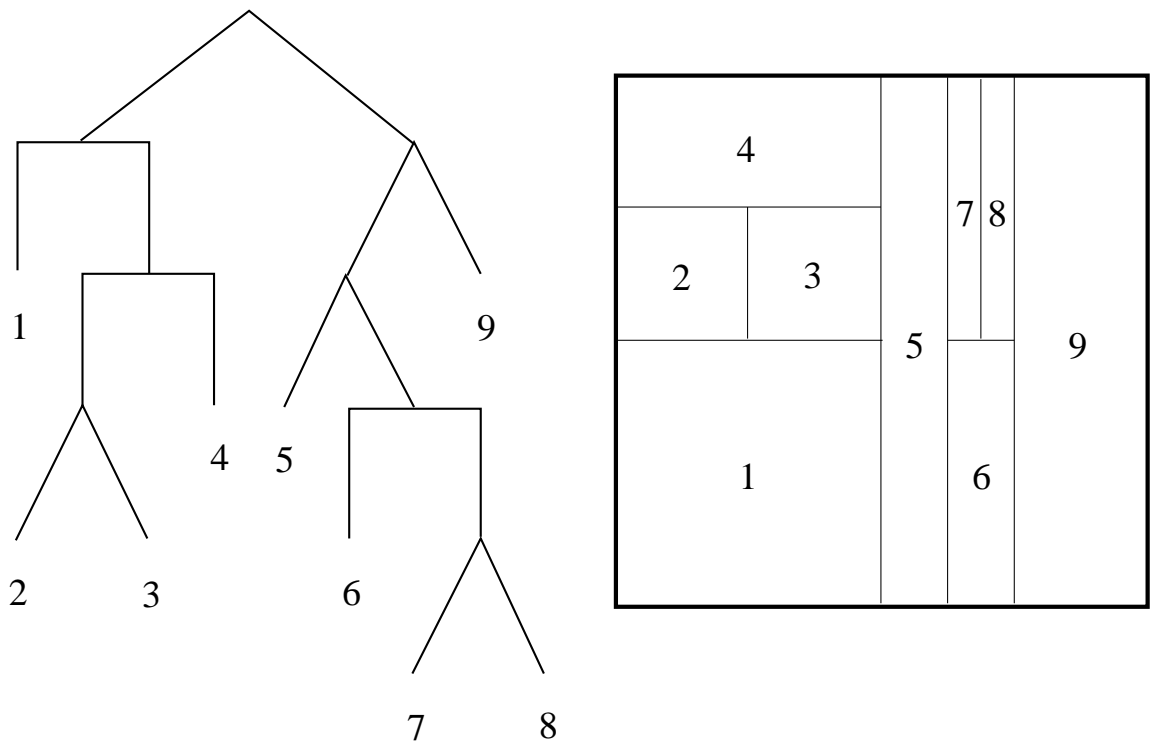}\\
\caption{An example of a binary tree with the two types of carets,
and its corresponding partition of $I^2$.} \label{tree-div}
\end{figure}

See Figure \ref{tree-div} for an example of a partition of $I^2$ and
its corresponding binary tree. Every dyadic partition of the unit
square can be obtained with a binary tree with these two types of
carets.

For the general groups $nV$, the partitions are divisions of the
unit $n$-cube in $n$-rectangles of dyadic lengths, and the binary
trees have $n$ types of carets. For simplicity, we state
our results for $2V$, but all of the results
below extend naturally to the cases for $nV$, with $n>1$.

\begin{figure}
\includegraphics[width=4in]{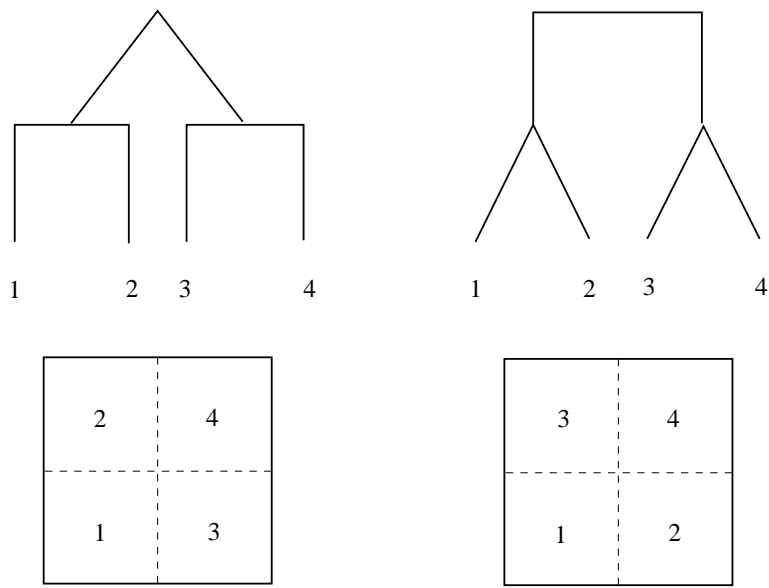}\\
\caption{The relation obtained when performing vertical and
horizontal subdivisions} \label{relation}
\end{figure}

The fact that we have vertical and horizontal subdivisions brings
new relations to the group. These relations arise when both types of
subdivisions are combined in different orders to obtain different
descriptions of the same dyadic partition of $I^2$. The obvious relation
 (and the one from which all other relations are deducible)
  is the combination of one subdivision of
each type in the two possible orders, as illustrated in Figure \ref{relation}.

\begin{figure}
\includegraphics[width=4in]{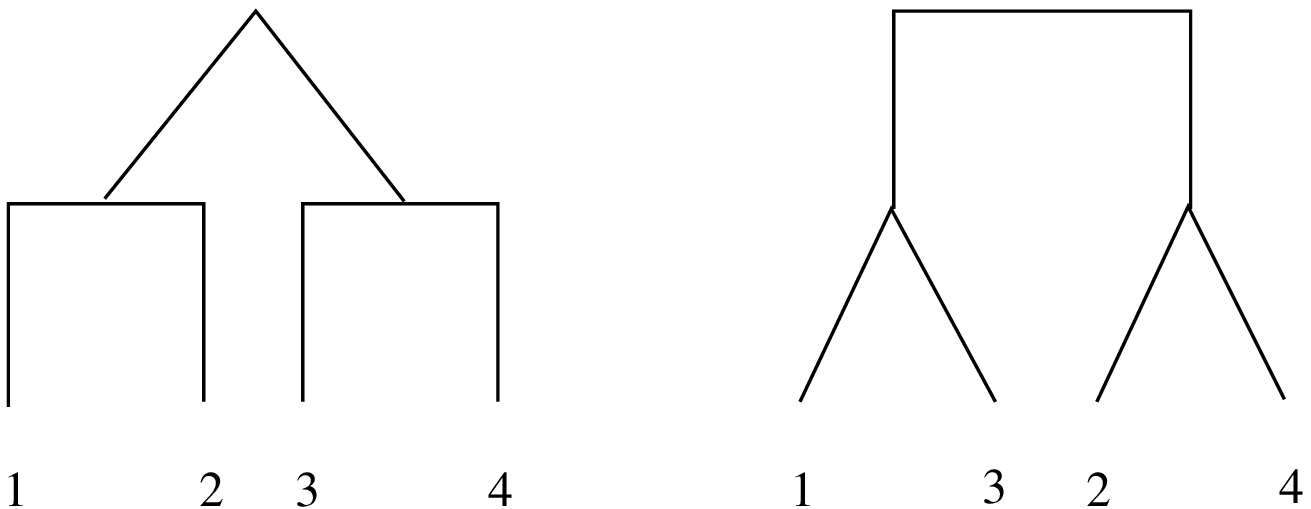}\\
\caption{A tree-pair diagram for the identity element which has a
non-identity permutation} \label{weird-id}
\end{figure}

Subdividing in both the  vertical and horizontal directions once,
but in the two possible orders, gives  the same dyadic
partitions, but
according to our convention the resulting rectangles are numbered in
different way. This fact makes using tree diagrams to express
the multiplication of elements tricky and quite unwieldy for large
elements.  Even though the leaves of the two-caret-type tree diagrams
are ordered in a natural way, this order is not apparent in the
square, and it is not preserved when different diagrams represent the
same partition. Thus, there are diagrams representing the identity
whose leaves are ordered in different ways--- an example is illustrated
in  Figure
\ref{weird-id}.

\section{Families of generators}

Brin \cite{brinhigher1,brinhigher2} showed that these
groups $nV$ are finitely generated and  gave several
generating sets for $2V$. As it is common with groups of the
Thompson family, there is an infinite presentation, which is useful for its
symmetry and regularity, which has a finite subpresentation.
For the purposes of word length, we work with the standard
finite generating set described by Brin.
We will next define the families of generators for $2V$. As it is
customary in Thompson's groups, see \cite{cfp}, the generators are built on a
backbone of an all-right tree with only vertical carets.
\begin{enumerate}
\item The generators $A_n$ involve only vertical subdivisions, and
they are the traditional generators of Thompson's group $F$.

\item The generators $B_n$ have the relevant caret replaced by a
horizontal one.

\item The generators $C_n$ have the last caret of the all-right tree
of horizontal type, and they are used to build horizontal carets on
the right-hand-side of the tree, as will be seen later.

\item The generators $\pi_n$ and $\overline\pi_n$ are permutations
built on an all-right tree with only vertical carets. These
generators are exactly equal to those for the subgroup $V$ appearing as the
purely vertical elements.
\end{enumerate}

\begin{figure}
\includegraphics[width=3.2in]{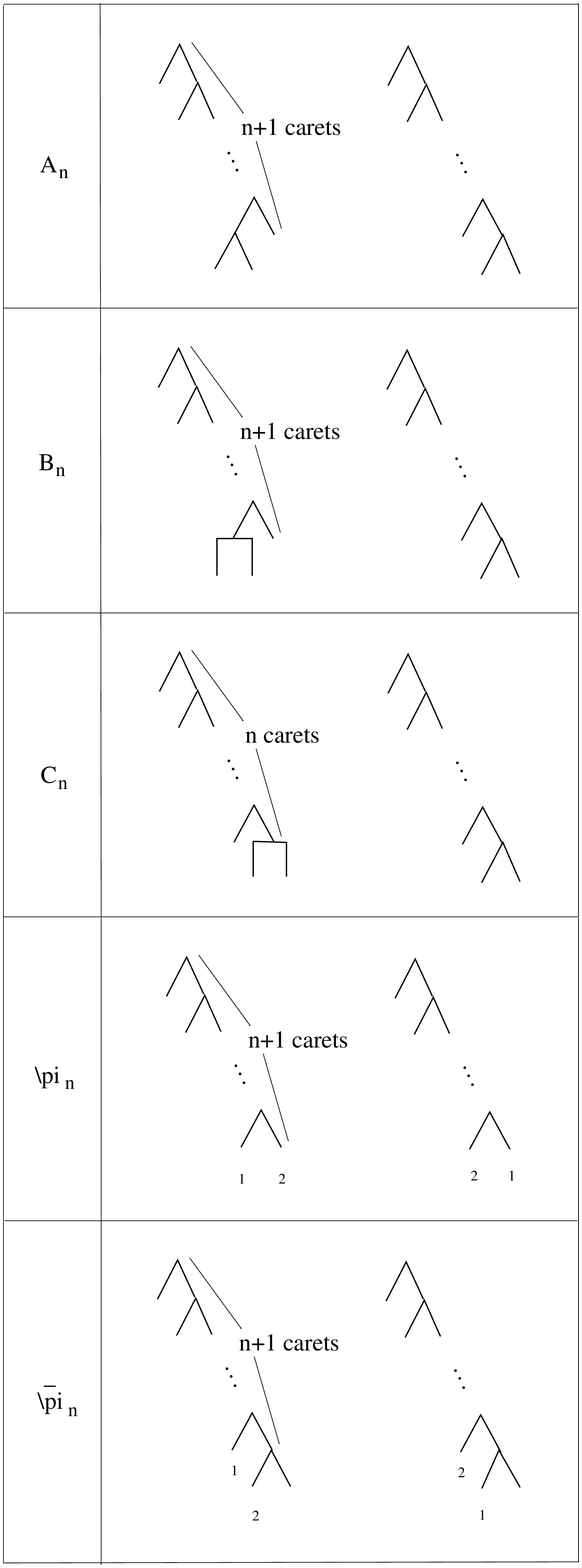}\\
\caption{Generators of $2V$} \label{gens}
\end{figure}

\begin{thm}[Brin \cite{brinhigher1,brinhigher2}] The families
$A_n$, $B_n$, $C_n$, $\pi_n$ and $\overline\pi_n$ generate $2V$.
\end{thm}

{\it Proof:} Since our proof is actually quite different from that
of Brin's, we will include it here, and use aspects of it later for metric considerations.

An element of $2V$ is now a triple $(T_+,\pi,T_-$), where the two
trees are composed of the two types of carets. First, we  subdivide a given element
into three elements: $(T_+,id,R_k)$, $(R_k,\pi,R_k)$, and
$(R_k,id,T_-)$, where $R_k$ is the all-right tree with $k$ leaves and only vertical carets.
Clearly $(R_k,\pi,R_k)$ is product of the permutation generators, as
it is the case in $V$ already.

To obtain the element $(T_+,id,R_k)$, we will concentrate first on
the backbone of the tree $T_+$; that is, the sequence of carets in
the right-hand-side. If this sequence of carets has horizontal
carets in the positions $m_1,m_2,\ldots,m_p$, then the product of
the generators $C_{m_1}C_{m_2}\ldots C_{m_p}$ produces exactly a
backbone with horizontal carets in the desired positions.

\begin{figure}
\includegraphics[width=5in]{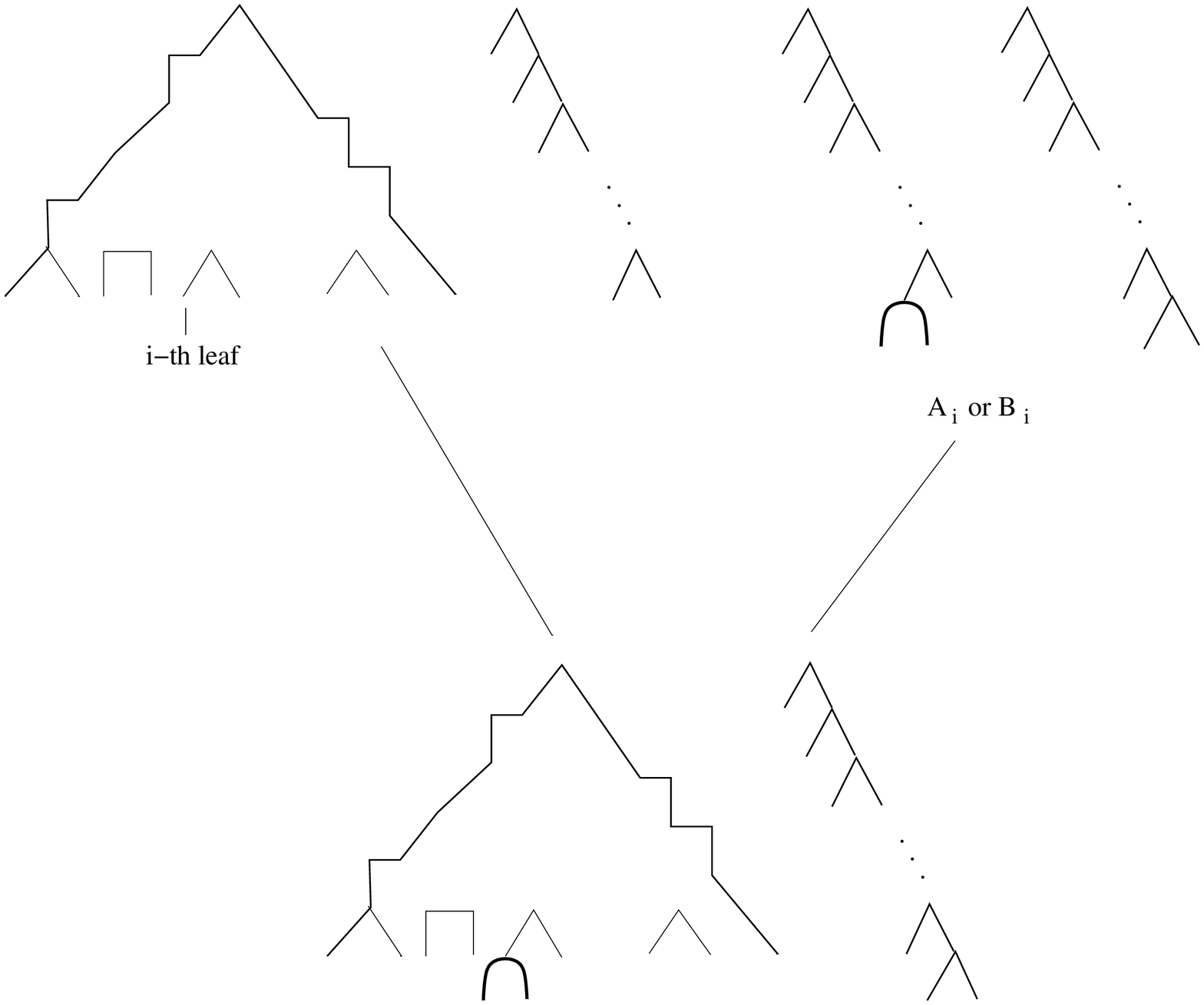}\\
\caption{The building process: attaching a caret to the $i$-th leaf
by multiplying by $A_i$ or $B_i$.} \label{multi}
\end{figure}

\begin{figure}
\includegraphics[width=6in]{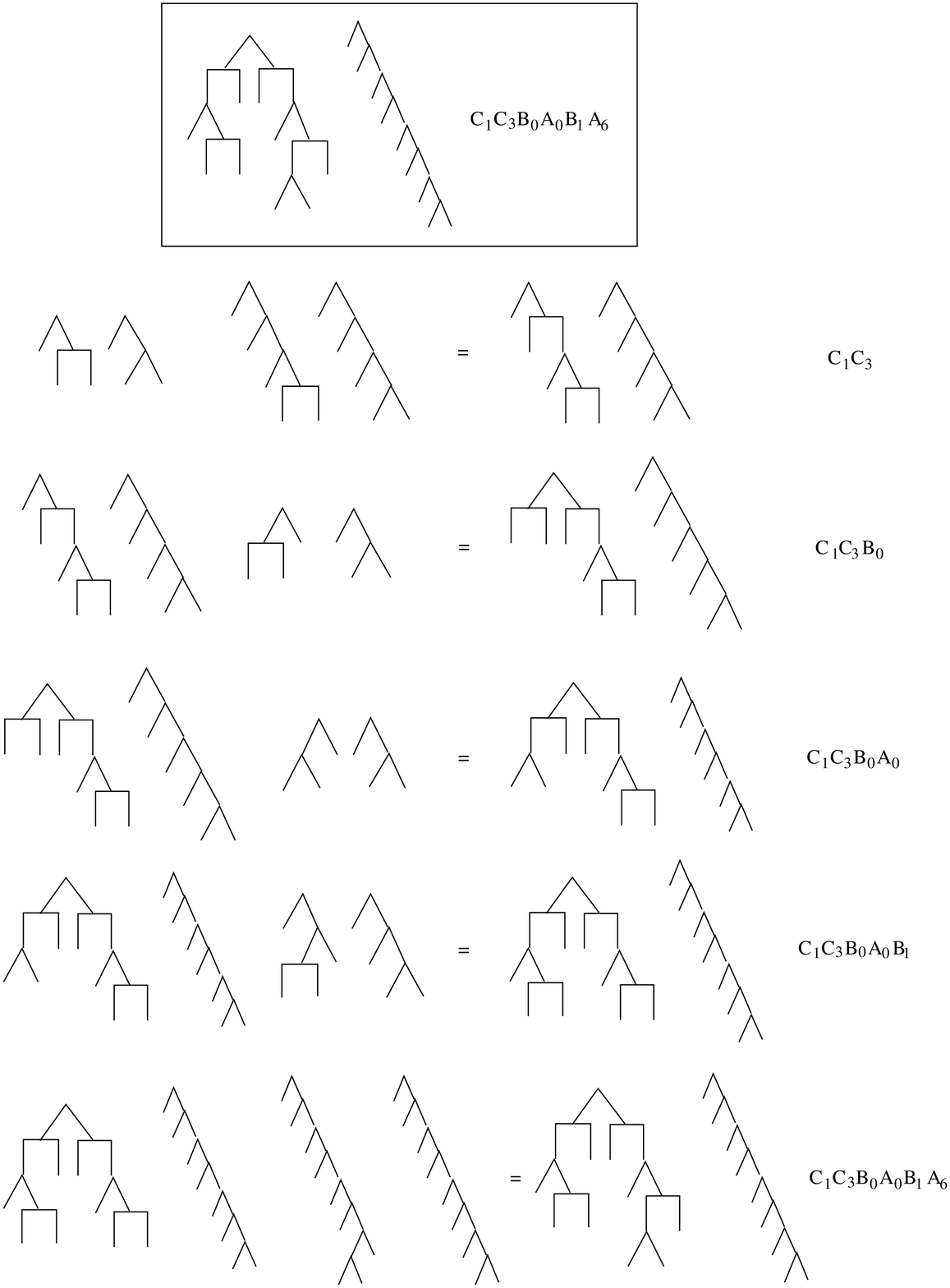}\\
\caption{An example of construction of an element generator by
generator.} \label{word}
\end{figure}

Once the backbone is constructed, then each new caret is obtained by
a generator of the type $A_i$ or $B_i$. To attach a vertical caret
from the leaf labelled $i$, we only need to multiply by $A_i$.
Similarly, to attach a  horizontal caret to leaf $i$, we multiply by $B_i$. This
 multiplication process is shown in Figure \ref{multi}, with an example
 illustrated in
Figure \ref{word}.

This proves that the element $(T_+,id,R_k)$ is product of the
generators $A_i$, $B_i$ and $C_i$, and using inverses, we have that
the entire group is generated by the full family of generators. \qed

An element of the type $(T_+,id,R_k)$ is called a \emph{positive} element of $2V$.
Positive elements can always be written as products of the generators $A_i$, $B_i$ and $C_i$ without using their inverses.

In the process of proving Brin's theorem, we have obtained the
following result:

\begin{thm}\label{semi-normal form}
For an element of $2V$, with respect to the standard infinite generating set
$\{A_i, B_i, C_i, \pi_i, \overline\pi_n\}$, we have:
\begin{itemize}
\item A positive element of $2V$ always admits an expression of
the type
$$
C_{m_1}\ldots C_{m_p} W_1(A_{i_1},B_{i_1})\ldots W_r(A_{i_r},B_{i_r}).
$$
where the $W_i$ are words on the positive generators $A_i$ and $B_i$ and never their inverses, and also where
$$
m_1<m_2<\ldots m_p\qquad
i_1<i_2<\ldots i_r
$$
\item Any element of $2V$ always admits an expression $P\Pi Q^{-1}$, where $P$ and $Q$ are positive elements, and $\Pi$ is a permutation on an all-right vertical tree, and thus a word in the $\pi_i$ and $ \overline\pi_i$.
\end{itemize}
\end{thm}

This expression will be used as a semi-normal form for elements of
$2V$.

\section{Metric properties}

We will be interested of metric properties up to the standard affine equivalence,
defined here.

\begin{defn} Given two functions
$$
f,g:G\longrightarrow\mathbb{R},
$$
we say that $f\prec g$ if there exists a constant $C>0$ such that for all $x\in G$, we have that $f(x)\le Cg(x   )$. Also, we say that $f$ and $g$ are equivalent, written $f\sim g$ if $f\prec g$ and $g\prec f$.
\end{defn}

Elements of Thompson's groups admit representations as diagrams with binary trees in one form or another, perhaps with  associated permutations or braids between the leaves. For several of these groups, the distance $|x|$ of an element $x$ to the identity is closely related with the number of carets $N(x)$ of the minimal diagram--- indicating that the  complexity of the minimal (reduced) diagram is a good indication of the size of the element.

In Thompson's groups $F$ and $T$, the two functions are equivalent, for both groups we have that $|x|\sim N(x)$, see \cite{bcs} and \cite{thompt}. This result cannot hold in $V$, merely for
counting reasons. The number of elements of $V$ whose trees have $N$ carets is at least the order of $N!$ (due to the permutations), while the number of distinct elements of length $N$ can be at most exponential in $N$ in any group.

The best possible bounds for the number of carets in terms of word length in  the group $V$ is the following inequality proved by Birget as Theorem 3.8  \cite{birget}. We have
$$
N(x) \prec |x| \prec N(x) \log N(x),
$$
and these bounds are optimal, in a sense similar to that which will be made precise in Section \ref{bounds}.

Since the lower bounds on word length are linear in the number of carets in the case of $V$  we see that the standard inclusions $F\subset T\subset V$ are all quasi-isometric embeddings; that is, $F$ is undistorted in $T$, and $T$ is undistorted in $V$. This property will be no longer true for inclusions in $2V$, as we will show in Section \ref{distortion}.

In this section we will prove the analog for $2V$ to the inequality above:

\begin{thm}\label{logn-nlogn} In the group $2V$, the word length of an element $|x|$ and number of carets $N(x)$ in a minimal size tree-pair representative satisfy
$$
\log N(x) \prec |x| \prec N(x) \log N(x),
$$
and the bounds cannot be improved.
\end{thm}

{\it Proof:} The upper bound is proven the standard way. We take a positive word $P$, and according to Proposition \ref{semi-normal form}, write it as
$$
C_{m_1}\ldots C_{m_p} W_1(A_{i_1},B_{i_1})\ldots W_r(A_{i_r},B_{i_r}).
$$
Then, we rewrite each generator using the following relations:
$$
A_{i+1}=A_0^{-i+1}A_1A_0^{i-1}\qquad B_{i+1}=A_0^{-i+1}B_1A_0^{i-1}\qquad C_{i+1}=A_0^{-i+1}C_1A_0^{i-1}
$$
for all $i>1$. Since the conjugating element is always $A_0$, cancellations ensure that the length stays approximately the same. The bound $N\log N$ appears because of the permutations, since $N\log N$ is the diameter of $\mathcal S_N$ with respect to the relevant transpositions.

\begin{figure}
\includegraphics[width=3in]{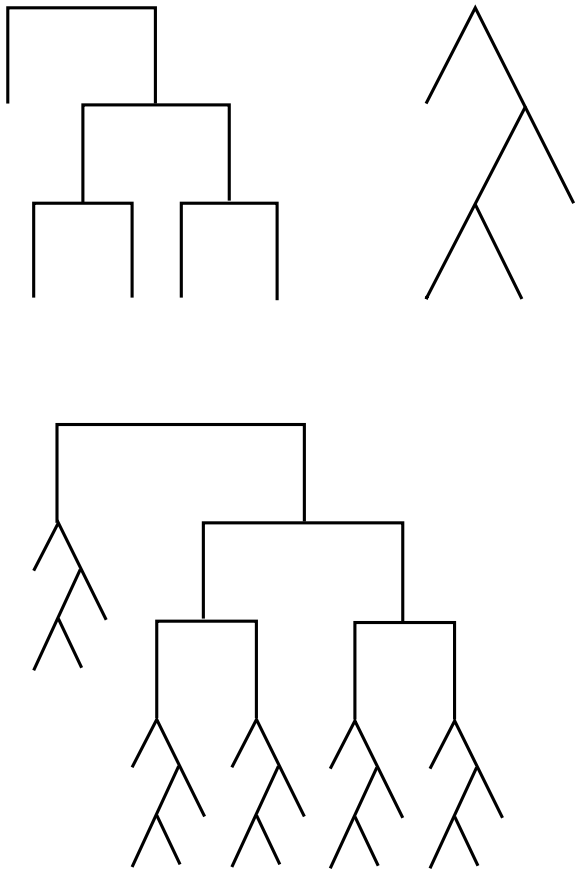}\\
\caption{An illustration of why when two types of carets are involved, the number of carets gets multiplied. The tree on the top left has only horizontal subdivisions, so when joined with a tree with only vertical ones, these have to be put at every leaf.} \label{log-bound}
\end{figure}

For the lower bound, we note that if an element has $k$ carets, when it is multiplied by a generator, it is possible that the multiplication could have up to $4k$ carets. If an element has, for instance, only horizontal subdivisions, when multiplied by $A_1$, the number of subdivisions (and hence the number of carets) becomes $4k$, since each one of the previous subdivision is divided in four, as illustrated in Figure \ref{log-bound}

If each multiplication by a generator can multiply the number of carets by 4, an element of length $\ell$ could have up to $4^{\ell}$ carets.  Other generators may increase word length by additive factors, but the
worst case is exactly the increase by a factor of 4.  Hence $N(x) \le 4^{\ell}$, so $\log N(x) \le  |x|$.\qed

\section{Optimality of bounds}\label{bounds}

First, we describe some properties of the lower bound for word length in terms of the
number of carets.

\begin{lem} The word length of an element of $2V$ which is represented by a tree-pair diagram of depth $n$ is at least $n/3$
\end{lem}

{\it Proof:}  Right multiplying by a generator in the finite generating set can add most 3 levels to the tree-pair diagram.
This happens in  the case that we right-multiply an element $(T_+, \pi, T_-)$ with no right subtree in $T_-$ by any of $A_1, B_1$ or $C_1$ from
the finite generating set.  Thus,
any product of less than $n/3$ generators cannot have depth $n$. \qed

From this it follows that the cyclic subgroups generated by any single generator
from the $A_n$ or $B_n$ families are
undistorted in $2V$, for example.

As seen in the proof of Theorem \ref{logn-nlogn}, multiplying by a generator can
increase the number of carets by a multiplicative factor. The powers of
the element $C_0$ have  exponentially
many carets, in fact the minimal number of carets to represent $C_0^n$ is $2^n$.  We give further
related specific examples later in Section \ref{distortion}. Any of these examples show that
 the lower bound of Theorem \ref{logn-nlogn} is optimal.

To analyze the genericity of the upper bound in Theorem \ref{logn-nlogn}, we use arguments analogous to those
for $V$.
Birget \cite{birget} showed not only the $n \log n$ upper bound on growth word length in $V$ with respect to the number of carets, but also that the fraction of elements close to this bound converges
exponentially fast to 1.  In 2V, we note analogous behavior:

\begin{thm}
Let $H_n$ be the set of elements of $2V$ representable with $n$ carets and
having no representatives with fewer than $n$ carets.  The
fraction of elements of $H_n$ which have word length greater than $n \log n$ converges
exponentially fast to 1.
\end{thm}

{\it Proof:}
Here we use a counting argument, analogous to that used by Birget \cite{birget} to count elements
in $V$.  In $2V$, we consider the subset of  representatives of elements which have only horizontal subdivisions in the domain and only vertical subdivisions in the range.  Their tree-pair diagrams are guaranteed to be reduced and of minimal size by a straightforward argument.  So to count the
set of elements of this type of diagram size $n$, we have $C_n$ choices (where $C_n$ is the $n$-th Catalan number) for the first all-square tree of size $n$ and $C_n$ choices for the all-triangular tree of size $n$, and $n!$ choices for the permutation.  This set of elements has size $C_n^2 n!$ which is, by the Stirling formula:

$$(C_n)^2 n! = \left( \frac{(2n)!}{(n+1)! n!}\right)^2 n!  = \sqrt{\frac2{\pi}} \frac{16^ n}{e^n} n^{(n-2)} (1+o(1))
$$

where the $o(1)$ term goes to zero as $n$ increases.

 The number of elements of word length $n$ in any finitely generated group with respect to $d$ generators is no more than $(2d)^n$.  Thus we see that the ratio of elements that have word length less than $n \log_{2d} n$ out of elements that have tree-pair diagrams of size $n$ is less than $$ \frac{n^n}{  \sqrt{\frac2{\pi}} \frac{16^ n}{e^n} n^{(n-2)} (1+o(1))} \sim c \frac{n^2 e^n}{16^n}$$
 which converges to 0 exponentially fast, as desired.  So by complementing we have the result.

\section{Distortion}\label{distortion}

Elements of $2V$ can be represented with pairs of binary trees, where the carets have been subdivided in two types. This fact makes it more difficult to multiply elements whose caret types disagree, illustrated already in Theorem \ref{logn-nlogn}, and the number of carets can grow faster because of this situation. This phenomenon has been described already by Wladis in \cite{wladisdiss} for the group $F(2,3)$, which has also carets of two types (binary and ternary).  This feature of $2V$ implies now that the groups $F$, $T$ and $V$, when seen as subgroups of $2V$ (by doing only vertical subdivisions, for instance), are exponentially distorted.

\begin{thm} The groups $F$, $T$ and $V$ are at least exponentially distorted in $2V$. \end{thm}

\begin{figure}
\includegraphics[width=6in]{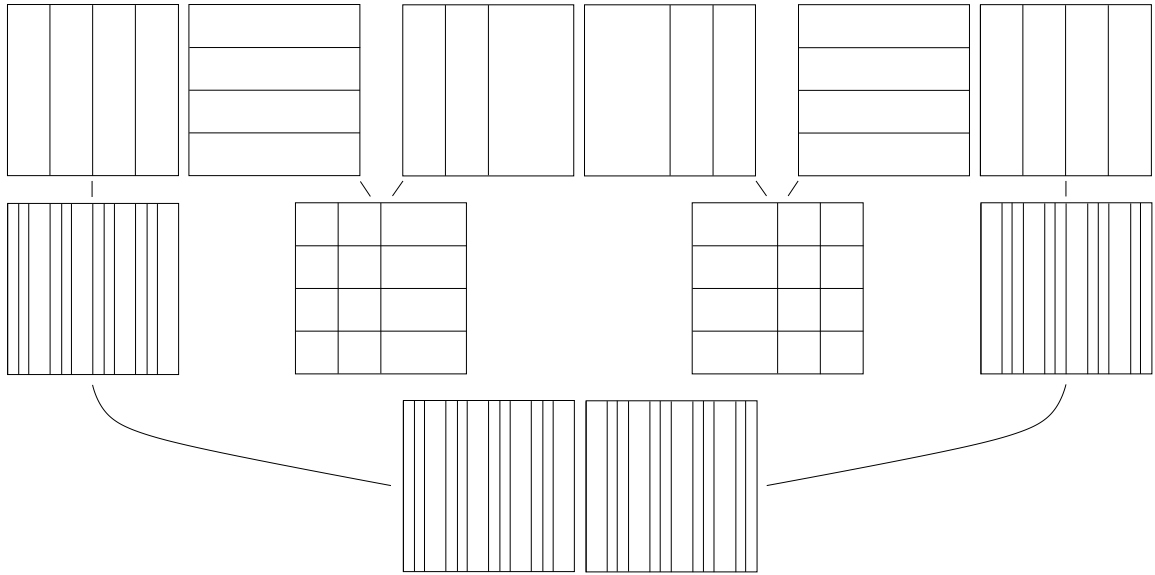}\\
\caption{The process of multiplication for $C_0^{-2}X_0C_0^2$, illustrating the exponential number of subdivisions.} \label{expdist}
\end{figure}

{\it Proof.}
We consider the specific element $C_0^{-n}X_0C_0^n$, illustrated in Figure  \ref{expdist}.
This element lies in a copy of  $F$ in $2V$ obtained by putting $F$ into $2V$ using
only vertical subdivisions.

The element $C_0^n$ has $2^n$ carets, as it is seen with an easy induction. Its two trees are balanced trees of depth $n$, one with only vertical subdivisions, and one with only horizontal subdivisions. By matching the horizontal subdivisions with the vertical ones in $X_0$ we obtain that the element $C_0^{-n}X_0C_0^n$ has a number of carets of the order of $2^n$, and all of them of vertical type, so the element is in $F$.  The element $C_0^{-n}X_0C_0^n$ has length in $V$ no more than $2n+1$,
but the number of carets is exponential in $N$ and thus its word length as an element of the
vertical copy of $F$ in $2V$ is also exponential.

Thus $F$ is at least exponentially distorted in $2V$. Since $F$ is undistorted in $T$ and $V$,
then $T$ and $V$ are also at least exponentially distorted in $2V$. \qed

\bibliographystyle{plain}

\end{document}